\documentclass[12pt,a4paper,reqno]{amsart}
\usepackage[english]{babel}
\usepackage{amssymb,latexsym,amsfonts,amsthm,upref,amsmath}
\usepackage[foot]{amsaddr}
\usepackage[margin=1in]{geometry} 
\usepackage[dvipsnames,x11names]{xcolor}
 \definecolor{myblue}{HTML}{003399}
\usepackage{hyperref}
\hypersetup{colorlinks,citecolor=myblue,filecolor=black,linkcolor=myblue,urlcolor=myblue}
\usepackage{enumerate}  
\usepackage{tikz}
\usepackage{float}
\usepackage{cleveref}
\usepackage{mathtools}
\usepackage{cite}
\usepackage{filecontents}
\usepackage{enumitem}
\makeatletter
\newcommand{\leqnomode}{\tagsleft@true}
\newcommand{\reqnomode}{\tagsleft@false}

\makeatother
\newtheorem*{thm*}{Theorem}
\newtheorem*{lem*}{Lemma}
\newtheoremstyle{prim}{}{}{\normalfont}{}{\bfseries}{.}{ }{}
\newtheoremstyle{stil}{}{}{\slshape}{}{\bfseries}{.}{ }{}
\theoremstyle{stil}
\newtheorem{thm}{Theorem}[section]
\newtheoremstyle{defi}{}{}{}{}{\bfseries}{.}{ }{}
\theoremstyle{defi}

\theoremstyle{defi}

\theoremstyle{stil}
\newtheorem*{mthm*}{Main Theorem}
\newtheorem*{kor*}{Corollary}
\newtheorem{pro}[thm]{Proposition}
\theoremstyle{stil}
\newtheorem{lem}[thm]{Lemma}
\theoremstyle{stil}

\theoremstyle{prim}

\newenvironment{prf}{\noindent \textit{Proof.}}{\null\hfill$\qed$\hskip
2mm\vskip 2mm}
\DeclareMathOperator{\A}{A}
\newcommand{\dygc}{  {\rm DY} (\mathfrak{gl}_{m|n})}
\newcommand{\dyg}{ {\rm DY}(\mathfrak{gl}_{m|n})}
\newcommand{\dys}{ {\rm DY}(\mathfrak{sl}_{n|n})}
\newcommand{\dyp}{ {\rm DY}(\A(n-1,n-1))}
\newcommand{\Y}{ {\rm Y}(\mathfrak{gl}_{m|n})}
\newcommand{\Yd}{ {\rm Y}^+(\mathfrak{gl}_{m|n})}

\newcommand{\U}{ {\rm U}}
\newcommand{\glmn}{\mathfrak{gl}_{m|n}}
\newcommand{\slnn}{\mathfrak{sl}_{n|n}}
\newcommand{\CC}{\mathbb{C}}
\newcommand{\ZZ}{\mathbb{Z}}
\newcommand{\wvr}{\overline}

\newcommand{\ot}{\otimes}
\newcommand{\ts}{\hspace{1pt}}

\newcommand{\sgn}{ \mathop{\rm sgn}}

\newcommand{\diag}{\mathop{\mathrm{diag}}}

\newcommand{\fand}{\quad\text{and}\quad}
\newcommand{\Fand}{\qquad\text{and}\qquad}

\newcommand{\non}{\nonumber}
\newcommand{\beq}{\begin{equation}}
\newcommand{\eeq}{\end{equation}}
\newcommand{\ben}{\begin{equation*}}
\newcommand{\een}{\end{equation*}}

\makeatletter
\def\smalloverbrace#1{\mathop{\vbox{\m@th\ialign{##\crcr\noalign{\kern3\p@}%
  \tiny\downbracefill\crcr\noalign{\kern3\p@\nointerlineskip}%
  $\hfil\displaystyle{#1}\hfil$\crcr}}}\limits}
\makeatother

\makeatletter
\def\smallunderbrace#1{\mathop{\vtop{\m@th\ialign{##\crcr
   $\hfil\displaystyle{#1}\hfil$\crcr
   \noalign{\kern3\p@\nointerlineskip}%
   \tiny\upbracefill\crcr\noalign{\kern3\p@}}}}\limits}
\makeatother

\begin{document}

\title{A note on the quantum Berezinian for the double Yangian  of the Lie superalgebra  $\mathfrak{gl}_{m|n}$}

\author{Lucia Bagnoli}
\author{Slaven Ko\v{z}i\'{c}}
\address[L. Bagnoli and S. Ko\v{z}i\'{c}]{Department of Mathematics, Faculty of Science, University of Zagreb,  Bijeni\v{c}ka cesta 30, 10\,000 Zagreb, Croatia}
\email{lucia.bagnoli@math.hr}
\email{kslaven@math.hr}

 \begin{abstract}
In this note, we generalize the notion of quantum Berezinian to the double Yangian ${\rm DY}(\mathfrak{gl}_{m|n})$  of the Lie superalgebra  $ \mathfrak{gl}_{m|n} $. We show that its coefficients form a   family of algebraically independent topological generators of the center of ${\rm DY}(\mathfrak{gl}_{m|n})$.
\end{abstract}

\maketitle

\allowdisplaybreaks

\section{Introduction}\label{intro}
\numberwithin{equation}{section}

The  {\em Yangian} $\Y$ for the  general linear Lie superalgebra $\mathfrak{gl}_{m|n}$ is a deformation of the universal enveloping algebra $\U(\glmn[t])$. It was introduced by Nazarov,   who also proved that it possesses a certain remarkable family of central elements; see \cite{Naz}.  They arise as  coefficients of the  formal power series $b_{m|n}^- (u)$, called {\em quantum Berezinian}, which  is a super analogue of the  quantum determinant   for the  Yangian of $  \mathfrak{gl}_{n} $; cf. \cite[Sect. 2]{MNO}. It was conjectured by Nazarov \cite{Naz}  and   proved by Gow \cite{Gow} that the  coefficients of $b_{m|n}^- (u)$ generate the    center of  $\Y$. Later on, different aspects of the quantum Berezinian and its applications were extensively studied; see, e.g., the papers by Lu \cite{Lu2209}, Lu and Mukhin \cite{LM}, Molev and Ragoucy \cite{MR}, Tsymbaliuk \cite{T} and references therein. Moreover, a generalization of the quantum Berezinian to arbitrary parities was introduced and investigated; see the papers by  Chang and Hu \cite{CH} and   Huang and Mukhin \cite{HM}.

In this note,   we consider the   {\em double Yangian} $\dyg$ for  the    Lie superalgebra $\mathfrak{gl}_{m|n}$, defined via an $R$-matrix presentation, which goes back to Zhang \cite{Z}, and our goal is to determine its center. The double Yangian
contains two important subalgebras, the    Yangian  $\Y$ and the {\em dual Yangian} $\Yd$. Motivated by \cite{Naz}, we introduce the quantum Berezinian $b_{m|n}^+ (u)$ for the dual Yangian and then we prove that the coefficients of  $b_{m|n}^- (u)$ and $b_{m|n}^+ (u)$ are algebraically independent topological generators of the center of $\dyg$. Furthermore, we show that, for $m=n$, this result naturally extends to the double Yangian $\dys$ for the  special linear Lie superalgebra $\mathfrak{sl}_{n|n}$.

To establish these results, we  follow the approach of Gow \cite{Gow0,Gow}, who found a new proof of aforementioned Nazarov's theorem  \cite{Naz} on centrality of the coefficients of $b_{m|n}^- (u)$.
It employs, in particular, super analogues of certain techniques and results on the Yangian for $ \mathfrak{gl}_{n} $,  going back to the papers of  Brundan and Kleshchev \cite{BK}, Molev, Nazarov and Olshanski \cite{MNO} and Nazarov and Tarasov \cite{NT}, which we adapt to our setting.
 One key step in this approach is construction of the Drinfeld generators of the double Yangian, as the proof takes place in the corresponding realization. As with the cases of Yangian   \cite{Gow} and the  double Yangian for    $ \mathfrak{gl}_{n} $   \cite{I,JY}, the generators are found via   Gauss decomposition which relies on the theory of quasideterminants \cite{GGRW,GR}. Finally, the fact that the  coefficients of the series $b_{m|n}^- (u)$ and $b_{m|n}^+ (u)$ generate the entire center of $\dyg$   is verified using the Poincar\'{e}--Birkhoff--Witt Theorem  for the double Yangian \cite{BaKo}.

\section{Preliminaries on the double Yangian for \texorpdfstring{$\glmn$}{glm|n}}\label{section1}
 
Consider the loop Lie superalgebra  $\mathcal{L}(\glmn)=\glmn\ot\CC[t,t^{-1}]$   of $\glmn$. 
Its supercommutation relations are given by
$$
[e_{ij}(r),e_{kl}(s)]
= 
\delta_{kj}\ts  e_{il}(r+s)
-  \delta_{il}\ts e_{kj}(r+s) (-1)^{(\bar{i}+\bar{j})(\bar{k}+\bar{l})}, 
$$
where $e_{ij}(r)=e_{ij}\ot t^r$ and
$e_{ij}\in \glmn$ are  matrix units. 
The element
$e_{ij}(r)$ is of parity $\bar{i}+\bar{j}$, where  $\bar{i}=0$ for $i=1,\ldots ,m$ and $\bar{i}=1$ for $i=m+1,\ldots ,m+n$.
To express the defining relations for the double Yangian, we shall need the rational $R$-matrix 
\beq\label{Rmatrix}
R(u)=1-hPu^{-1},\quad\text{where}\quad 
1=\sum_{i,j=1}^{m+n} e_{ii}\ot e_{jj}
\fand
P=\sum_{i,j=1}^{m+n} e_{ij}\ot e_{ji}\ts (-1)^{\bar{j}}.
\eeq

The
({\em level $0$}) {\em double Yangian $\dyg$ for $\mathfrak{gl}_{m|n}$} is defined as the $\ZZ_2$-graded unital associative algebra over the commutative ring $\CC[[h]]$ generated by the elements $t_{ij}^{(\pm r)}$ of parity $\bar{i}+\bar{j}$, where $  i,j=1,\ldots , m+n$ and $r=1,2,\ldots,$ subject to the  defining relations which are written in terms of the generator matrices 
$$
T^{-}(u)=\sum_{i,j=1}^{m+n} (-1)^{\bar{i}\bar{j}+\bar{j}} e_{ij}\ot t^{-}_{ij}(u)
\Fand
T^+ (u)=\sum_{i,j=1}^{m+n} (-1)^{\bar{i}\bar{j}+\bar{j}} e_{ij}\ot t^+_{ij}(u).
$$
Their matrix entries are given by
\beq\label{series}
t_{ij}^-(u)=\delta_{ij}+h\sum_{r\geqslant 1} t_{ij}^{(r)} u^{-r}
\Fand
t_{ij}^+ (u)=\delta_{ij}-h\sum_{r\geqslant 1} t_{ij}^{(-r)} u^{r-1} .
\eeq
Finally, the defining relations take the form
\begin{align}
&R(u-v)\ts T^{\pm}_1(u)\ts T^{\pm}_2(v)
=T^{\pm}_2(v)\ts T^{\pm}_1(u)\ts R(u-v)\label{RTT},\\
&R(u-v)\ts T^-_1 (u)\ts T_2^+ (v)
=T_2^+ (v)\ts T^-_1 (u)\ts R(u-v). \label{DYRTT}
\end{align}

Denote by $ \Y$ (resp. $\Yd$)  the {\em Yangian} (resp. the {\em dual Yangian}), i.e. the  unital   subalgebra of the double Yangian generated by   all elements $t_{ij}^{(r)}$ (resp. $t_{ij}^{(-r)}$) with $i,j=1,\ldots ,m+n$ and $r=1,2,\ldots .$ Throughout the paper, we assume that all these algebras are complete with respect to the $h$-adic topology. 
 Observe that the matrices $T^\pm (u)$ are invertible, which in the case of $T^+(u)$ is due to the fact that $\dyg$ is $h$-adically completed.
Let us write
$$
T^\pm (u)^{-1}=\sum_{i,j=1}^{m+n} (-1)^{\bar{i}\bar{j}+\bar{j}} e_{ij}\ot t^{\pm\prime}_{ij} (u).
$$
By employing the $RTT$-relations \eqref{RTT} and \eqref{DYRTT}, one   derives  the formulae
	\begin{align}
	&(u-v)[t^\pm_{ij} (u),t^{\pm\prime}_{kl}  (v)]
= (-1)^{\bar{i}\bar{j}+\bar{i}\bar{k}+\bar{j}\bar{k} }h
\left(\hspace{-2pt}\delta_{kj}\sum_{a=1}^{m+n} t^\pm_{ia}(u)\ts t^{\pm\prime}_{al} (v)
-\delta_{il}\sum_{a=1}^{m+n} t^{\pm\prime}_{ka} (v)\ts t^\pm_{aj}(u)\hspace{-2pt}\right),\label{comm1}\\
&(u-v)[t^\pm_{ij} (u),t^{\mp\prime}_{kl} (v)]
= (-1)^{\bar{i}\bar{j}+\bar{i}\bar{k}+\bar{j}\bar{k}} h
\left(\hspace{-2pt}\delta_{kj}\sum_{a=1}^{m+n} t^\pm_{ia}(u)\ts t^{\mp\prime}_{al} (v)
-\delta_{il}\sum_{a=1}^{m+n} t^{\mp\prime}_{ka} (v)\ts t^\pm_{aj}(u)\hspace{-2pt}\right).\label{comm2}
	\end{align}

By setting $h=1$ in \eqref{Rmatrix}--\eqref{DYRTT}, one obtains the definition of the double Yangian for $\mathfrak{gl}_{m|n}$  over the complex field. Define the degrees of its generators by 
$
\deg t_{ij}^{(r)}=r-1$ and $\deg t_{ij}^{(-r)}=-r$ 
for all $i,j=1,\ldots ,m+n$ and $r=1,2,\ldots .$ 
By the Poincar\'{e}--Birkhoff--Witt Theorem for the double Yangian for $\mathfrak{gl}_{m|n}$ over $\CC$ \cite[Thm. 2.6]{BaKo}, its corresponding graded algebra is isomorphic to the universal enveloping algebra ${\rm U}(\mathcal{L}(\glmn))$. The isomorphism  is given by the    assignments
$$\bar{t}_{ij}^{(r)}\mapsto (-1)^{\bar{i}}\ts e_{ij}(r-1)\fand \bar{t}_{ij}^{(-r)}\mapsto (-1)^{\bar{i}}\ts e_{ij}(-r)$$
 with $r\geqslant 1$, where $\bar{t}_{ij}^{(\pm r)}$ indicate the images of the generators $t_{ij}^{(\pm r)}$ in the respective components of the corresponding graded algebra. It 
 naturally extends to the case of double Yangian $\dygc$ over $\CC[[h]]$ via classical limit $h\to 0$, so that we have
\beq\label{isomorphism}
\dygc / h\dygc \cong  {\rm U}(\mathcal{L}(\glmn)) . 
\eeq
Indeed, this is due to the fact that the top degree components with respect to the degree operator $\deg$ in the defining relations  over $\CC$   coincide with the coefficients of the lowest power of the parameter $h$ in the corresponding relations over $\CC[[h]]$.\footnote{We remark that, in particular,  \cite[Thm. 2.6]{BaKo} holds for $m=n$, so that we have the same isomorphism \eqref{isomorphism} in this case. Indeed, throughout the entire paper \cite{BaKo}, it is assumed  that $m\neq n$, as its main results   are proved only in that case. However, \cite[Thm. 2.6]{BaKo} holds 
for $m=n$ as well and this can be verified by repeating the arguments from its proof in the $m\neq n$ case, which are presented therein.}

 \section{The center of the double Yangian for \texorpdfstring{$\glmn$}{glm|n}}\label{section2}
In this section, we prove our main result on the center of the double Yangian $\dyg$, Theorem \ref{centpro} below. In Subsection \ref{sub21}, we present the Gauss decomposition for the double Yangian, which we need in Subsection \ref{sub22} to obtain its Drinfeld generators; see Theorem \ref{lem_34}. In Subsection \ref{sub23}, we study the quantum Berezinian and, finally, prove our main result. At the end, in Subsection \ref{sub24}, we extend Theorem \ref{centpro} to the case of  double Yangian for $\slnn$.

\subsection{Gauss decomposition}\label{sub21} 

 The matrices $T^{\pm}(u)$ possess the {\em Gauss decomposition}
\beq\label{gaussdec}
T^{\pm}(u)=F^{\pm}(u)\ts D^{\pm}(u)\ts E^{\pm }(u),
\eeq
where $F^{\pm}(u)$ and $E^{\pm}(u)$ are  triangular matrices with units on their main diagonals,
\begin{align*}
F^{\pm}(u)=\begin{pmatrix}
1&0&\ldots &0\\
f^{\pm}_{21}(u)&1&\ldots&0\\
\vdots&\vdots&\ddots&\vdots\\
f^{\pm}_{m+n\ts 1}(u)&f^{\pm}_{m+n\ts 2}(u)&\ldots&1
\end{pmatrix},
\qquad
E^{\pm}(u)=\begin{pmatrix}
1&e^{\pm}_{12}(u)&\ldots &e^{\pm}_{1\ts m+n}(u)\\
0&1&\ldots&e^{\pm}_{2\ts m+n}(u)\\
\vdots&\vdots&\ddots&\vdots\\
0&0&\ldots&1
\end{pmatrix},
\end{align*}
and $D^{\pm}(u)=\diag(d_1^{\pm}(u),\ldots, d_{m+n}^{\pm}(u))$ are diagonal matrices. The existence of  decomposition \eqref{gaussdec}, which  goes back to  \cite{Gow0,Gow, Z},
 is a consequence of the well-known results   from the theory of quasideterminants   \cite[Thms. 4.9.6, 4.9.7]{GGRW} (see also \cite{GR}), which also imply the explicit  expressions for the    entries of the above matrices as follows (see also \cite[Chaps. 1.10, 1.11]{M_book}). First, let $A=(a_{ij})$ be an arbitrary  $n\times n$ matrix over a unital ring.
Denote by $A^{ij}$ its submatrix obtained by deleting its $i$-th row and $j$-th column and suppose that $A^{ij}$ is invertible. One defines  the {\em $(i,j)$-th quasideterminant}   by 
\renewcommand*{\arraystretch}{0.6}
\setlength\arraycolsep{2pt}
$$\left|A\right|_{ij}=
\left|\begin{matrix}
a_{11}&\ldots &a_{1j}&\ldots &a_{1n}\\
\vdots& &\vdots& &\vdots\\
a_{i1}&\ldots &\boxed{a_{ij}}&\ldots &a_{in}\\
\vdots& &\vdots& &\vdots\\
a_{n1}&\ldots &a_{nj}&\ldots &a_{nn}
\end{matrix}\right|
\coloneqq
a_{ij} - r_i^j(A^{ij})^{-1}c_j^i,
$$
where $r_i^j$ (resp. $c_j^i$) is the row matrix (resp. column matrix) obtained from the $i$-th row (resp. $j$-th column) of $A$ by deleting the element $a_{ij}$.
The matrix entries of $D^{\pm}(u)$, $F^{\pm}(u)$ and $E^{\pm }(u)$ are found by
\renewcommand*{\arraystretch}{1.0}
\setlength\arraycolsep{5pt}
\begin{align}
&d_i^\pm(u)=\left|\begin{matrix}
t^\pm_{11}(u)&\ldots&t^\pm_{1\ts i-1}(u)& t^\pm_{1\ts i}(u)\\
\vdots& &\vdots&\vdots\\
t^\pm_{i\ts 1}(u)&\ldots&t^\pm_{i\ts i-1}(u)& \boxed{t^\pm_{i\ts i}(u)}
\end{matrix}\right|,\label{de}\\
&f_{ji}^\pm(u)=\left|\begin{matrix}
t^\pm_{11}(u)&\ldots&t^\pm_{1\ts i-1}(u)& t^\pm_{1\ts i}(u)\\
\vdots& &\vdots&\vdots\\
t^\pm_{i-1\ts 1}(u)&\ldots&t^\pm_{i-1\ts i-1}(u)&  t^\pm_{i-1\ts i}(u) \\
t^\pm_{j\ts 1}(u)&\ldots&t^\pm_{j\ts i-1}(u)&  \boxed{t^\pm_{j\ts i}(u)}
\end{matrix}\right|d_i^\pm (u)^{-1},\label{ef}\\
&e_{ij}^\pm(u)=d_i^\pm (u)^{-1}\left|\begin{matrix}
t^\pm_{11}(u)&\ldots&t^\pm_{1\ts i-1}(u)& t^\pm_{1\ts j}(u)\\
\vdots& &\vdots&\vdots\\
t^\pm_{i-1\ts 1}(u)&\ldots&t^\pm_{i-1\ts i-1}(u)&  t^\pm_{i-1\ts j}(u) \\
t^\pm_{i\ts 1}(u)&\ldots&t^\pm_{i\ts i-1}(u)&  \boxed{t^\pm_{i\ts j}(u)}
\end{matrix}\right|.\label{e}
\end{align} 
One easily checks that the series in \eqref{de}--\eqref{e} are of the form
\begin{align*}
&f_{ij}^- (u)=h\sum_{r\geqslant 1} f_{ij}^{(r)}u^{-r},
&&e_{ij}^- (u)=h\sum_{r\geqslant 1} e_{ij}^{(r)}u^{-r},
&&d_i^- (u)=1+h\sum_{r\geqslant 1} d_i^{(r)}u^{-r},
\\
&f_{ij}^+ (u)=-h\sum_{r\geqslant 1} f_{ij}^{(-r)}u^{r-1},
&&e_{ij}^+ (u)=-h\sum_{r\geqslant 1} e_{ij}^{(-r)}u^{r-1},
&&d_i^+ (u)=1-h\sum_{r\geqslant 1} d_i^{(-r)}u^{r-1},
\end{align*}
where the coefficients $f_{ij}^{(\pm r)},e_{ij}^{(\pm r)}, d_i^{(\pm r)}$ belong to the double Yangian.
In particular, for $i=1,\ldots ,m+n-1$ and $r=1,2,\ldots ,$ we  write
$f_{i }^{(\pm r)}=f_{i+1\ts i}^{(\pm r)}$, $e_{i }^{(\pm r)}=e_{i\ts i+1 }^{(\pm r)}$,
\beq\label{generatori_2}
f_{i}^\pm (u)=f_{i+1\ts i}^\pm (u)\fand 
e_i^\pm (u)=e_{i\ts i+1}^\pm (u).
\eeq
We point out that from \eqref{gaussdec} it is easy to recover the following relations, for $1 \leqslant i<j \leqslant m+n $, (cf. \cite[Section 4]{Gow}) 
\begin{align}
\begin{split}
\label{equazioni_t}
&t_{ii}^\pm(u)=d^\pm_{i}(u)+\sum_{k<i}f^\pm_{ik}(u)d^\pm_{k}(u)e^\pm_{ki}(u),\\
&t_{ij}^\pm(u)=d^\pm_{i}(u)e^\pm_{ij}(u)+\sum_{k<i}f^\pm_{ik}(u)d^\pm_{k}(u)e^\pm_{kj}(u),\\
&t_{ji}^\pm(u)=f^\pm_{ji}(u)d^\pm_{i}(u)+\sum_{k<i}f^\pm_{jk}(u)d^\pm_{k}(u)e^\pm_{ki}(u),\\
&t^{\pm\prime}_{ii}(u)=d^\pm_{i}(u)^{-1}+\sum_{k>i}e^{\pm\prime}_{ik}(u)d^\pm_{k}(u)^{-1}f^{\pm\prime}_{ki}(u),\\
&t^{\pm\prime}_{ij}(u)=e^{\pm\prime}_{ij}(u)d^\pm_{j}(u)^{-1}+\sum_{k>j}e^{\pm\prime}_{ik}(u)d^\pm_{k}(u)^{-1}f^{\pm\prime}_{kj}(u),\\
&t^{\pm\prime}_{ji}(u)=d^\pm_{j}(u)^{-1}f^{\pm\prime}_{ji}(u)+\sum_{k>j}e^{\pm\prime}_{jk}(u)d^\pm_{k}(u)^{-1}f^{\pm\prime}_{ki}(u),
\end{split}
\end{align}
where the series $e^{\pm\prime}_{ij}(u)$ and $f^{\pm\prime}_{ji}(u)$  are given by
\begin{align*}
&e^{\pm\prime}_{ij}(u)=\sum_{i=i_0<i_1<...<i_{s}=j}(-1)^{s}e^{\pm }_{i_0 i_1}(u)e^{\pm }_{i_1 i_2}(u) \cdots e^{\pm }_{i_{s-1} i_s}(u),\\
&f^{\pm\prime}_{ji}(u)=\sum_{i=i_0<i_1<...<i_{s}=j}(-1)^{s}f^{\pm }_{i_s i_{s-1}}(u) \cdots f^{\pm }_{i_2 i_1}(u) f^{\pm }_{i_1 i_0}(u).
\end{align*}
\subsection{Drinfeld generators}\label{sub22} 
In the next three lemmas, we study certain maps   which will be instrumental in the proofs of the main results of this paper.
First, by using   defining relations \eqref{RTT} and \eqref{DYRTT},    one easily verifies the next lemma, which extends the maps from \cite[Lemma 4.1]{Gow} and \cite[Prop. 4.2]{Gow} to the double Yangian.

\begin{lem}\label{lem_auto}
There exist    unique associative algebra isomorphisms 
$$\rho_{m|n}\colon\dygc\to {\rm DY} (\mathfrak{gl}_{n|m})\Fand \omega_{m|n}\colon\dygc\to\dygc$$  
such that for all $i,j=1,\ldots ,m+n$ we have
\begin{align*}
&\rho_{m|n}(t_{ij}^\pm (u))
=t_{m+n+1-i\ts\ts m+n+1-j}^\pm (-u)\Fand \omega_{m|n}(T^\pm (u))
=T^\pm (-u)^{-1}.
\end{align*}
\end{lem}

Consider the isomorphism $\zeta_{m|n}\colon\dygc\to  {\rm DY} (\mathfrak{gl}_{n|m})$ defined by 
$$\zeta_{m|n}= \rho_{m|n}\circ \omega_{m|n}. $$

\begin{lem}\label{lem_31}
For all    $i=1,\ldots ,m+n-1$ and $j=1,\ldots ,m+n$ we have
\begin{gather*}
\zeta_{m|n}(f_i^\pm (u))=-e_{m+n-i}^\pm (u) ,\quad
\zeta_{m|n}(e_i^\pm (u))=-f_{m+n-i}^\pm (u) ,\\
\zeta_{m|n}(d_j^\pm (u))=d_{m+n-j+1}^\pm (u)^{-1}.
\end{gather*}
\end{lem}

\begin{prf}
The statement of  lemma for the series $f_i^- (u)$, $e_i^- (u)$ and $d_j^-(u)$  is due to \cite[Prop. 4.2]{Gow}. Its proof  directly generalizes to the   remaining series $f_i^+ (u)$, $e_i^+ (u)$ and $d_j^+(u)$.
\end{prf}

For any positive integer $k$ let $\iota_{m|n,k}\colon {\rm DY} (\mathfrak{gl}_{m|n})\hookrightarrow  {\rm DY} (\mathfrak{gl}_{m+k|n})$  be the inclusion, i.e.
\beq\label{iota}
\iota_{m|n,k}(t_{ij}^{(\pm r)})=t_{i+k\ts j+k}^{(\pm r)}\quad\text{for all }i,j=1,\ldots ,m+n,\,r=1,2,\ldots.
\eeq
Consider the algebra monomorphism
$\psi_{m|n,k}\colon  {\rm DY} (\mathfrak{gl}_{m|n})\to {\rm DY} (\mathfrak{gl}_{m+k|n})$, defined  by 
\beq\label{psi}
\psi_{m|n,k}=\omega_{m+k|n}\circ \iota_{m|n,k}\circ \omega_{m|n},
\eeq
which generalizes the  monomorphism  ${\rm Y}(\mathfrak{gl}_{m})\to  {\rm Y}(\mathfrak{gl}_{m+k})$
 of Nazarov and Tarasov  \cite{NT}.

\begin{lem}\label{lem_32}
For all    $i=1,\ldots ,m+n-1$, $j=1,\ldots ,m+n$ and $k=1,2,\ldots$ we have
\begin{gather*}
\psi_{m|n,k}(f_i^\pm (u))=f_{k+i}^\pm (u) ,\quad
\psi_{m|n,k}(e_i^\pm (u))=e_{k+i}^\pm (u) ,\quad
\psi_{m|n,k}(d_j^\pm (u))=d_{k+j}^\pm (u).
\end{gather*}
\end{lem}

\begin{prf}
The lemma is a consequence of the formula for the image of $t_{pq}^\pm (u)$ under   $\psi_{m|n,k}$,
\beq\label{psiformula}
\psi_{m|n,k}(t_{pq}^\pm (u))
=\left|\begin{matrix}
t_{11}^\pm(u)&\ldots&t_{1k}^\pm(u)&t_{1\ts k+q}^\pm(u)\\
\vdots& &\vdots&\vdots\\
t_{k1}^\pm(u)&\ldots&t_{kk}^\pm(u)&t_{k\ts k+q}^\pm(u)\\
t_{k+p\ts 1}^\pm(u)&\ldots&t_{k+p\ts k}^\pm(u)&\boxed{t_{k+p\ts k+q}^\pm(u)}\\
\end{matrix}\right|   \text{ for all }p,q=1,\ldots ,m+n.
\eeq
The even analogue of the above expression for the image of $ t_{pq}^- (u) $ was found by Brundan and Kleshchev  \cite[Lemma 4.2]{BK}. Their argument naturally extends to the super case, both for  $ t_{pq}^- (u) $ and $ t_{pq}^+ (u) $; see also \cite[Sect. 4]{Gow}.
\end{prf}

We are now in position to establish the Drinfeld generators for the double Yangian.

\begin{thm}\label{lem_34}
The coefficients
\beq\label{coefficients}
f_i^{(\pm r)},\,e_i^{(\pm r)},\, d_j^{(\pm r)},\quad\text{where}\quad i=1,\ldots ,m+n-1,\, j=1,\ldots ,m+n,\, r=1,2,\ldots,
\eeq
topologically generate the   double Yangian $\dygc$.
\end{thm}

\begin{prf}
It is known by \cite[Sect. 3]{Gow} that the coefficients of $f_i^-(u)$, $e_i^-(u)$ and $d_j^- (u)$ generate the Yangian $\Y$, so it is sufficient to check that the remaining coefficients    generate the dual Yangian $\Yd$. In other words, it suffices to show that the coefficients of all $f_{pq}^+(u)$ and $e_{qp}^+(u)$
with $1\leqslant q<p\leqslant m+n$ can be expressed in terms  of \eqref{coefficients}. We shall present the proof for  $e_{qp}^+(u)$. The case of $f_{pq}^+(u)$ follows by suitably adjusting the same arguments. 
The proof goes by induction over $m+n$. Clearly, for $m+n=2$, the thesis holds.
Suppose the thesis holds for the coefficients of all $e_{qp}^+(u)$ and $f_{pq}^+(u)$ in $\dygc$ such that $m+n\leqslant t$ for some integer $t \geqslant 2$. Let us consider $\dygc$ such that $m+n=t+1$. Using \eqref{e} and \eqref{psiformula} one   shows that 
\beq
\label{psimap2}
\psi_{m-1|n,1}(e_{ij}^+(u))=e_{i+1\ts j+1}^+(u)\qquad\text{for all}\quad 1\leqslant i<j\leqslant m+n-1. 
\eeq
Then by induction the thesis holds for all $e_{i\ts j}^+(u)$ in $\dygc$ with $2\leqslant i<j\leqslant m+n$. Since ${\rm DY} (\mathfrak{gl}_{m|n-1}) \subset {\rm DY} (\mathfrak{gl}_{m|n})$, by induction the thesis holds for all $e_{i\ts j}^+(u)$ in $\dygc$ with $1\leqslant i<j\leqslant m+n-1$. Finally, by  using relation   \eqref{comm2} for $(i,j,k,l)=(1,t,t,t+1)$ and relations \eqref{equazioni_t}, one can express the coefficients of $e_{1\ts t+1}^+(u)$ in terms of  \eqref{coefficients}.   
\end{prf}

The next lemma is another key result   proved by using  properties of  $\psi_{m|n,k}$ and $\zeta_{m|n}$. Essentially, as with Proposition \ref{prop37} below, it follows by an argument from \cite[Sect. 4]{Gow}, which can be extended to other relations; see, e.g., \cite[Cor. 3.4]{Lu2308}. Nonetheless,  we provide detailed proofs for completeness.

\begin{lem}\label{lem_dovi}
For all $i,j=1,\ldots ,m+n$ we have
\beq\label{dscomm}
d_i^\pm(u)\ts d_j^{\pm}(v) =d_j^{\pm}(v)\ts d_i^\pm(u)
\Fand
d_i^-(u)\ts d_j^{+}(v) =d_j^{+}(v)\ts d_i^-(u).
\eeq
\end{lem}

\begin{prf}
Let us prove the first identity in \eqref{dscomm}. It is well-known that the quasideterminants $d_i^\pm(u)$ and $d_j^{\pm}(v)$ commute for $1\leqslant i,j\leqslant m$;  see, e.g., \cite[Thm. 3.3]{I} or \cite[Thm. 2.5]{JY}.
Furthermore, if
$m+1\leqslant i,j\leqslant m+n$, one concludes that $d_i^\pm(u)$ and $d_j^{\pm}(v)$ commute by applying the map $\psi_{0|n,m}$ on the identity
$
d_{i-m}^\pm(u)\ts d_{j-m}^{\pm}(v) =d_{j-m}^{\pm}(v)\ts d_{i-m}^\pm(u),
$
which holds in $ {\rm DY} (\mathfrak{gl}_{0|n})\cong  {\rm DY} (\mathfrak{gl}_{n})$, and then using Lemma \ref{lem_32}. Finally, suppose that $1\leqslant i <m+1\leqslant  j\leqslant m+n$. By \eqref{de}, the coefficients of $d_i^\pm (u)$ can be expressed in terms of  
\beq\label{temp1}
t_{pq}^{(\mp r)},\qquad\text{where}\quad 1\leqslant p,q\leqslant m\fand r=1,2,\ldots.
\eeq
On the other hand, by Lemma \ref{lem_32}, we have  $d_j^\pm (v)=\psi_{0|n,m}(d_{j-m}^\pm (v))$. Thus, we see from \eqref{iota} and \eqref{psi} that the coefficients of $d_j^\pm (v)$ can be expressed in terms of  
\beq\label{temp2}
t_{m+p\ts m+q}^{(\mp r)\prime} ,\qquad\text{where}\quad  1\leqslant p,q\leqslant  n\fand r=1,2,\ldots.
\eeq
As the elements in \eqref{temp1} and \eqref{temp2} commute by \eqref{comm1}, so do $d_i^\pm(u)$ and $d_j^{\pm}(v)$.

Let us prove the second identity in \eqref{dscomm}. Again, it is well-known that $d_i^-(u)$ and  $d_j^{+}(v)$ commute for $1\leqslant i,j\leqslant m$;  see \cite[Thm. 3.3]{I} or  \cite[Thm. 2.5]{JY}. If $1\leqslant i <m+1\leqslant  j\leqslant m+n$, we   conclude as before  that the coefficients of $d_i^-(u)$ can be expressed in terms of
$$
t_{pq}^{(  r)},\qquad\text{where}\quad 1\leqslant p,q\leqslant m\fand r=1,2,\ldots
$$
and the coefficients of $d_j^+(v)$ in terms of
\beq\label{telements}
t_{m+p\ts m+q}^{(- r)\prime} ,\qquad\text{where}\quad  1\leqslant p,q\leqslant  n\fand r=1,2,\ldots.
\eeq
Therefore,   $d_i^-(u)$ and  $d_j^{+}(v)$ commute due to \eqref{comm2}. It remains to verify the case $ m+1\leqslant i\leqslant m+n$, $1\leqslant j\leqslant m+n$. However, this immediately follows by applying the map $\zeta_{n|m}$ on the identities
$
d_{i-m}^-(u)\ts d_j^+(v)= d_j^+(v)\ts d_{i-m}^-(u)
$,
which hold in  $ {\rm DY} (\mathfrak{gl}_{n|m})$,
and then using Lemma \ref{lem_31}.
\end{prf}

\subsection{Quantum Berezinian}\label{sub23}  

Let $\mathfrak{S}_n$ be the symmetric group. Following \cite{Naz}, we define the {\em quantum Berezinian} for the matrices $T^\pm(u)$ by 
\begin{align*}
b_{m|n}^\pm (u)=&\sum_{\tau\in\mathfrak{S}_m }\sgn\tau\ts\ts t^\pm_{\tau(1)1}(u) \cdots t^\pm_{\tau(m)m}(u-(m-1)h  ) \\
&\times\sum_{\sigma\in\mathfrak{S}_n }\sgn\sigma\ts\ts t^{\pm\prime}_{m+1\ts m+\sigma(1)} (u-(m-1)h) \cdots t^{\pm\prime}_{m+n\ts m+\sigma(n)}(u-(m-n)h  ) .
\end{align*}
It is well-known (see, e.g., \cite{MNO,I}) that the first factors in the definition of $b_{m|n}^\pm (u)$,
\beq\label{quantumdeterminant}
C_m^\pm (u)=\sum_{\tau\in\mathfrak{S}_m }\sgn\tau\ts\ts t^\pm_{\tau(1)1}(u)t^\pm_{\tau(2)2}(u-h  )\cdots t^\pm_{\tau(m)m}(u-(m-1)h  ),  
\eeq
which coincide   with the images of the quantum determinants of $ {\rm DY} (\mathfrak{gl}_{m })\cong {\rm DY} (\mathfrak{gl}_{m|0})$ under the inclusion $ {\rm DY} (\mathfrak{gl}_{m|0})\hookrightarrow {\rm DY} (\mathfrak{gl}_{m|n})$,
can be written in the form
\beq\label{qdeterminant}
C_m^\pm (u)=d_1^\pm(u)\ts d_2^\pm(u-h)\cdots d_m^\pm(u-(m-1)h). 
\eeq
An analogue of this property for the quantum Berezinian is given by the next lemma.

\begin{lem}\label{pro_ber}
We have
\begin{align}
b_{m|n}^\pm (u)=&\,
d_1^\pm(u) \cdots d_m^\pm(u-(m-1)h)\non\\
&\times d_{m+1}^\pm(u-(m-1)h)^{-1} \cdots d_{m+n}^\pm(u-(m-n)h)^{-1}.\label{ber_d}
\end{align}
\end{lem}

\begin{prf}
The identity for $b_{m|n}^- (u)$ goes back to \cite[Thm. 1]{Gow0} and the one for $b_{m|n}^+ (u)$ can be verified by analogous arguments. They rely on the  expression \eqref{qdeterminant} for the quantum determinant and the observation that the quantum Berezinian can be written in the form
$$
b_{m|n}^\pm (u)=C_m^\pm (u)\ts\zeta_{n|m}(C_n^\pm(u-(m-n)h)),
$$
where $C_n^\pm(u )$ are   the images of the quantum determinants of $ {\rm DY} (\mathfrak{gl}_{n })\cong {\rm DY} (\mathfrak{gl}_{n|0})$ under the inclusion $ {\rm DY }(\mathfrak{gl}_{n|0})\hookrightarrow {\rm DY} (\mathfrak{gl}_{n|m})$; recall \eqref{quantumdeterminant}.
\end{prf}

The next proposition generalizes  Nazarov's theorem  on centrality of the coefficients of $b_{m|n}^- (u)$ \cite{Naz}  to the case of double Yangian. Its proof  follows the approach of Gow, who found a second    proof of this result in  \cite{Gow0}.

\begin{pro}\label{prop37}
The coefficients of $b_{m|n}^\pm (u)$ belong to the center of   $\dygc$.
\end{pro}

\begin{prf}
By Theorem \ref{lem_34}, it suffices to show that   the  series  in \eqref{de} and \eqref{generatori_2} commute with   $b_{m|n}^\pm (v)$. It is clear from Lemma \ref{lem_dovi} and the expression   for the quantum Berezinian  in \eqref{ber_d} that all $d_j^\pm (u)$ commute with both $b_{m|n}^+ (v)$ and $b_{m|n}^- (v)$, so  it is sufficient to consider the series $f_i^\pm (u)$ and  $e_i^\pm(u)$. 
Furthermore, due to \cite[Sect. 2]{Naz}, it is  known that all   $f_i^-(u)$ and  $e_i^-(u)$  commute with $b_{m|n}^- (v)$.  

\noindent(1) Let us prove that for $i=1,\ldots ,m+n-1$ we have
\beq\label{eq333}
e_i^+(u)\ts b_{m|n}^+ (v)=b_{m|n}^+ (v)\ts e_i^+(u)
\fand
f_i^+(u)\ts b_{m|n}^+ (v)=b_{m|n}^+ (v)\ts f_i^+(u).
\eeq

\noindent (1a) We start by verifying  \eqref{eq333} for $i=1,\ldots ,m-1$. Consider the   first equality. It is well-known from the even case that the quantum determinant \eqref{qdeterminant} commutes with $e_i^+(u)$ for $i=1,\ldots ,m-1$, as its coefficients belong to the center of  $ {\rm DY} (\mathfrak{gl}_{m })$. Hence, due to \eqref{ber_d}, it is sufficient to check the identities 
\beq\label{idds5}
e_i^+(u)\ts d_{m+j}^+(v)^{-1}=d_{m+j}^+(v)^{-1}e_i^+(u)\quad \text{for }j= 1,\ldots , n.
\eeq 
However, by the last family of equalities in Lemma \ref{lem_32}, we have
$\psi_{0|n,m}(d_{j}^+(v))=d_{m+j}^+(v)$. Hence, the coefficients of all $d_{m+j}^+(v)^{-1} $ with $j= 1,\ldots , n$ can be expressed in terms of elements \eqref{telements}. Therefore, the identities in \eqref{idds5} for $i=1,\ldots ,m-1$  follow immediately from \eqref{comm1}, so we conclude that the first equality in \eqref{eq333} for $i=1,\ldots ,m-1$ holds. The  second equality  can be proved analogously.

\noindent (1b) The equalities in \eqref{eq333} for $i= m+1,\ldots ,m+n-1$ are verified by arguing as in the corresponding part of the proof of \cite[Thm. 2]{Gow0}, so we omit the details.  However, it is worth noting that the argument relies on the properties of the map $\zeta_{m|n}$ from Lemma \ref{lem_31} and the   identities  in \eqref{eq333} for $i=1,\ldots ,m-1$.

\noindent (1c) Let $i=m$. We start by considering the case $m=n=1$ and prove that  
\beq\label{gl11}
e_1^+(u)\ts b_{1|1}^+ (v)=b_{1|1}^+ (v)\ts e_1^+(u)
\eeq
 holds in $ {\rm DY} (\mathfrak{gl}_{1|1})$. This is again checked in parallel with the corresponding part of the proof of \cite[Thm. 2]{Gow0}. More specifically, from Equality \eqref{comm1} for  $ {\rm Y}^+(\mathfrak{gl}_{1|1})$ with $(i,j,k,l)=(1,1,1,2)$ and  $(i,j,k,l)=(1,2,2,2)$ one derives   
$$
\left(u-v\right)e_1^+(v)\ts d_j^+(u)=\left(u-v-h\right)d_j^+(u)\ts e_1^+(v)+d_j^+(u)\ts e_1^+(u)\quad\text{for }j=1,2.
$$
By combining these two identities and the formula $b_{1|1}^+ (v)=d_1^+(v)d_2^+(v)^{-1}$, which follows from \eqref{ber_d}, one   proves \eqref{gl11}. 

Let us return to the general case of $\dygc$. By Lemma \ref{lem_32},   applying the map $\psi_{1|1,m-1}\colon  {\rm DY} (\mathfrak{gl}_{1|1})\to {\rm DY} (\mathfrak{gl}_{m|1})\subset  {\rm DY} (\mathfrak{gl}_{m|n})$ to \eqref{gl11} we obtain
\beq\label{glmny}
e_m^+(u)\ts d_m^+(v)\ts d_{m+1}^+(v)^{-1}=d_m^+(v)\ts d_{m+1}^+(v)^{-1}\ts e_m^+(u).
\eeq

Let us prove that $e_m^+(u)$ commutes with $d_j^+(v)$ for $j\neq m,m+1$. First, we  show that $e_1^+(u)$ commutes with $d_3^+(v),\ldots ,d_{n+1}^+(v)$ in $ {\rm DY} (\mathfrak{gl}_{1|n})$.  By \eqref{comm1},  in $ {\rm DY} (\mathfrak{gl}_{1|s})$  we have
\beq\label{eq21}
 t_{12}^+(u)\ts t_{s+1\ts s+1}^{+\prime}(v) =t_{s+1\ts s+1}^{+\prime}(v)\ts t_{12}^+(u)\quad\text{for all } s=2,\ldots ,n .
\eeq
As $t_{12}^+(u)=d_1^+(u)e_1^{+}(u)$ and $t_{s+1\ts s+1}^{+\prime}(v)=d_{s+1}^+(v)^{-1}$,
 by using Lemma \ref{lem_dovi}
one deduces from \eqref{eq21} that  in $ {\rm DY} (\mathfrak{gl}_{1|s})$ we have
\beq\label{eq22}
 e_1^+(u)\ts d_{s+1}^+(v)=d_{s+1}^+(v)\ts e_1^+(u) \quad\text{for all }s=2,\ldots , n.
\eeq
Note that the above equalities hold  in $ {\rm DY} (\mathfrak{gl}_{1|n})$ as well. Indeed, this follows by  applying the composition  $\zeta_{n|1}\circ\psi_{s|1,n-s}\circ \zeta_{1|s}$ to \eqref{eq22} and then employing Lemmas \ref{lem_31} and \ref{lem_32}.
Finally, by regarding \eqref{eq22} as identities in $ {\rm DY} (\mathfrak{gl}_{1|n})$ and applying the map $\psi_{1|n,m-1}$,  due to Lemma \ref{lem_32}, we immediately obtain the   identities in $\dygc$,
\beq\label{eq23}
 e_m^+(u)\ts d_{s+m}^+(v)=d_{s+m}^+(v)\ts e_m^+(u) \quad\text{for all }s=2,\ldots , n.
\eeq

By arguing as in the proofs of \eqref{glmny} and \eqref{eq23}, one can verify the equalities
\beq\label{eq33}
f_m^+(u)\ts d_m^+(v)\ts d_{m+1}^+(v)^{-1}=d_m^+(v)\ts d_{m+1}^+(v)^{-1}  f_m^+(u),
\quad
f_m^+(u)\ts d_{j}^+(v)=d_{j}^+(v)\ts f_m^+(u)
\eeq
in $\dygc$ with $j=m+2,\ldots ,m+n$.

Finally, in view of \eqref{glmny}, \eqref{eq23} and \eqref{eq33},
it remains to check that the series $e_m^+(u)$ and $ f_m^+(u)$ commute with $d_{1}^+(v),\ldots , d_{m-1}^+(v)$. However, this immediately follows by applying the map $\zeta_{n|m}$ to the identities
$$
e_{n}^+(u)\ts d_{s}^+(v)=d_{s}^+(v)\ts e_n^+(u)\fand
f_{n}^+(u)\ts d_{s}^+(v)=d_{s}^+(v)\ts f_n^+(u)\quad\text{with }s=n+2,\ldots, n+m,
$$
which, by the above discussion, hold in $ {\rm DY} (\mathfrak{gl}_{n|m})$, and then using Lemma \ref{lem_31}. Therefore, by the above arguments,  we conclude that   \eqref{eq333} holds for all  $i=1,\ldots ,m+n-1$.

\noindent(2) Our next goal is to prove that for $i=1,\ldots ,m+n-1$ we have
\beq\label{eq3333}
e_i^-(u)\ts b_{m|n}^+ (v)=b_{m|n}^+ (v)\ts e_i^-(u) 
\fand
f_i^-(u)\ts b_{m|n}^+ (v)=b_{m|n}^+ (v)\ts f_i^-(u).
\eeq

\noindent(2a) First, we consider the series $e_i^-(u)$ and $f_i^-(u)$ with $i=1,\ldots ,m-1$. The
equalities
\beq\label{qe1}
e_i^-(u)\ts d_j^+(v)=d_j^+(v)\ts e_i^-(u) \fand 
f_i^-(u)\ts d_j^+(v)=d_j^+(v)\ts f_i^-(u)
\eeq
for $j=m+1,\ldots ,m+n$
 follow  by an argument from the proof of Lemma \ref{lem_dovi}. More specifically, by \eqref{e}, the coefficients of  $e_i^-(u)$ and $f_i^-(u)$   can be expressed in terms of generators \eqref{temp1} while the coefficients of  $d^+_{m+1}(v),\ldots ,d^+_{m+n}(v)$  can be expressed in terms of elements \eqref{temp2}. Thus, they commute due to  \eqref{comm1}.
Next,  we observe that \eqref{qe1} holds for all $i=1,\ldots ,m-1$ and $j=1,\ldots ,m$ such that $j\neq i,i+1$ due to the presentation of the  double Yangian $ {\rm DY} (\mathfrak{gl}_{m})\cong {\rm DY} (\mathfrak{gl}_{m|0})\subset  {\rm DY} (\mathfrak{gl}_{m|n})$ \cite[Thm. 3.3]{I}; see also \cite[Thm. 2.5]{JY}. 
Finally, from the relations established in \cite[Thm. 3.3]{I}, one derives
$$
\left(u-v\right) \left(e_i^{-}(v)-e_i^+(v)\right)d_j^+(u)
=\left(u-v\mp h\right)d_j^+(u)\left(e_i^{-}(v)-e_i^+(v)\right)\quad\text{for }j=i,i+1,
$$
where the minus (resp. plus) sign  corresponds to the case $j=i$ (resp. $j=i+1$)\footnote{Observe that   the  notation of \cite[Thm. 3.3]{I} is connected with the setting of this paper via $X_i^+(v)\mapsto e_i^{-}(v)-e_i^+(v)$, $k_j^{-}(u)\mapsto d_j^+(u)$ and $\hbar\mapsto -h$.}. 
By employing the     above identities   along with 
$$
e_i^+(u)\ts d_i^+(v)\ts d_{i+1}^+(v-h)=d_i^+(v)\ts d_{i+1}^+(v-h)\ts  e_i^+(u),
$$
one can show that for all $i=1,\ldots ,m-1$  the series $e_i^-(u)$ commutes with the product $d_i^+(v)d_{i+1}^+(v-h)$. The same statement for $f_i^-(u)$  with $i=1,\ldots ,m-1$ is verified analogously.
  Hence, due to \eqref{ber_d}, the above discussion implies that \eqref{eq3333} holds for all $i=1,\ldots ,m-1$.
	
\noindent(2b) The fact that \eqref{eq3333} holds for $i=m+1,\ldots ,m+n$
follows  from the equalities in  \eqref{eq3333} for $i=1,\ldots ,m-1$
by the use of  properties of the isomorphism $\zeta_{m|n}$ given by Lemma \ref{lem_31}, which, in particular, imply, $\zeta_{m|n}(b_{m|n}^+ (u))=b_{n|m}^+ (u-(m-n)h)$.

\noindent(2c) It remains to check that  \eqref{eq3333} holds for $i=m$, but this can be done in parallel with the above proof of \eqref{eq333} in the case $i=m$.

\noindent(3) To prove the proposition, it remains to show that 
for $i=1,\ldots ,m+n-1$ we have
$$
e_i^+(u)\ts b_{m|n}^- (v)=b_{m|n}^- (v)\ts e_i^+(u) 
\fand
f_i^+(u)\ts b_{m|n}^- (v)=b_{m|n}^- (v)\ts f_i^+(u).
$$
However, this follows by suitably adapting the above proof of \eqref{eq3333}.
\end{prf}

To prove our main result, Theorem \ref{centpro} below, we shall need the next two lemmas. The first one is a direct generalization of \cite[Lemma 7.1]{BK} and it can be proved by  repeating the arguments of Brundan and Kleshchev; see also \cite[Prop. 2.12]{MNO}.

\begin{lem}\label{brkllema}
Let $\mathfrak{g}$ be a finite dimensional Lie algebra over $\CC$ and $\mathfrak{h}$ its reductive subalgebra. Denote by $\mathfrak{c}$   the centralizer of $\mathfrak{h}$ in $\mathfrak{g}$. The centralizer of the enveloping algebra $\U(\mathfrak{h}\ot\CC[t^{\pm 1}])$ in  $\U(\mathfrak{g}\ot\CC[t^{\pm 1}])$ is equal to $\U(\mathfrak{c}\ot\CC[t^{\pm 1}])$.
\end{lem}

The second lemma is verified by reducing the  problem to the even case of Lemma \ref{brkllema}. This can be  done in the same way as for $\glmn\ot\CC[t]$; see  \cite[Lemma 7.1]{Gow}.

\begin{lem}\label{brukle}
The center of ${\rm U}(\mathcal{L}(\glmn))$ is generated by the elements $I(r)= I\ot t^r$, $r\in\ZZ$.
\end{lem}

Let $b_{m|n}^{(\pm r)}$ for $r=1,2,\ldots$ be the coefficients of the quantum Berezinians so that
$$
b_{m|n}^{-}(u)=1+h\sum_{r\geqslant 1}b_{m|n}^{(  r)} u^{-r}
\fand
b_{m|n}^{+}(u)=1-h\sum_{r\geqslant 1}b_{m|n}^{(-r)} u^{r-1}.
$$

\begin{thm}\label{centpro}
The coefficients $b_{m|n}^{(\pm r)}$  with $r=1,2,\ldots$ are algebraically independent topological generators of the center of the double Yangian $\dygc$.
\end{thm}

\begin{prf}
The statement of the theorem  follows by the use of the isomorphism \eqref{isomorphism} and Lemmas \ref{brkllema} and \ref{brukle}. More precisely, we observe that the action of the isomorphism \eqref{isomorphism} on the classical limit of the coefficients of   quantum Berezinians is given by
\beq\label{hto0maps}
  b_{m|n}^{(  r)}  \big|_{h=0}\mapsto I(r-1)\fand
 b_{m|n}^{(  -r)} \big|_{h=0}\mapsto I(-r)
\eeq
for all $r=1,2,\ldots .$ Suppose $x$ is an element of the center of $\dyg$ such that $x|_{h=0}\neq 0$. Then the image of $x|_{h=0}$ under the map \eqref{isomorphism} is a nontrivial linear combination $\wvr{\lambda}_0$ of some monomials in $I(r)$. Let $\lambda_0$ be the corresponding linear combination of monomials in $ b_{m|n}^{(  \pm r)}$, such that the image of $\lambda_0|_{h=0}$ under \eqref{isomorphism} is $\wvr{\lambda}_0$. We have 
$$x-\lambda_0\in h \dyg,\quad\text{i.e.}\quad x=\lambda_0\mod h,$$
 and the element $x-\lambda_0$ again belongs to the center of the double Yangian. One can now proceed inductively as follows. The image of $h^{-1}(x-\lambda_0)|_{h=0}$ under   \eqref{isomorphism} is a  linear combination $\wvr{\lambda}_1$ of some monomials in $I(r)$. Let $\lambda_1$ be the corresponding linear combination of monomials in $ b_{m|n}^{(  \pm r)}$, such that the image of $\lambda_1|_{h=0}$ under \eqref{isomorphism} is $\wvr{\lambda}_1$. We have 
$$x-\lambda_0-h \lambda_1\in h^2 \dyg ,\quad\text{i.e.}\quad x=\lambda_0+h\lambda_1\mod h^2,$$ and the element $x-\lambda_0-h \lambda_1$ again belongs to the center of the double Yangian. With the double Yangian being topologically free, we can continue   the above procedure and express $x$ in terms of coefficients of quantum Berezinians $b_{m|n}^{\pm}(u)$.
\end{prf}

 \subsection{Application to the double Yangian for \texorpdfstring{$\slnn$}{sln|n}}\label{sub24} 

In this section, we establish an analogue of Theorem \ref{centpro} for the double Yangian of  $\slnn$. Introduce the series
\begin{align*}
&\delta_{i}^{-}(u)=1+h\sum_{r\geqslant 1} \delta_{i}^{(r)}u^{-r}=d_1^{-}(u)^{-1}d_{i+1}^{-}(u)\in {\rm DY} (\mathfrak{gl}_{n|n})[[u^{-1}]],\\
&\delta_{i}^{+}(u)=1-h\sum_{r\geqslant 1} \delta_{i}^{(-r)}u^{r-1}=d_1^{+}(u)^{-1}d_{i+1}^{+}(u)\in {\rm DY} (\mathfrak{gl}_{n|n})[[u]]
\end{align*}
for  $i=1,\ldots ,2n-1$. 
Motivated by \cite[Lemma 8.2]{Gow}, we define the {\em double Yangian $\dys$ for $\slnn$} as the $h$-adically completed subalgebra of ${\rm DY} (\mathfrak{gl}_{n|n})$ generated by the elements
\beq\label{generators_sl}
f_i^{(\pm r)},\,e_i^{(\pm r)},\, \delta_i^{(\pm r)},\quad\text{where}\quad i=1,\ldots ,2n-1, \, r=1,2,\ldots.
\eeq

 By examining the images of the generators \eqref{generators_sl} at $h\to 0$ under the isomorphism \eqref{isomorphism}, one finds that its restriction to   $\dys/h\dys$ produces the isomorphism
\beq\label{isomorphism_sl}
\dys / h\dys \cong  {\rm U}(\mathcal{L}(\slnn)) , 
\eeq
where $\mathcal{L}(\slnn)=\slnn\ot\CC[t,t^{-1}]$ is the loop Lie  superalgebra of $\slnn$. The coefficients of  quantum Berezinians belong to   $\dys$ as the series in \eqref{ber_d} can be expressed   as
\begin{align*}
b_{n|n}^{ \pm }(u)=&\,\delta^\pm _1(u-h)\ldots \delta^\pm_{n-1}(u-(n-1)h)\\
&\times \delta^\pm_{n}(u-(n-1)h)^{-1}\ldots \delta^\pm_{2n-2}(u-h)^{-1}\delta^\pm_{2n-1}(u)^{-1}.
\end{align*}

We recall from \cite{K} that the classical Lie superalgebra $\A(n-1,n-1)$ is defined by $\slnn/\mathbb{C}I$.
Let us denote by $\mathcal{L}(\A(n-1,n-1))=\A(n-1,n-1)\ot\CC[t,t^{-1}]$   its loop Lie  superalgebra. The next lemma is proved by the same arguments as Lemma \ref{brukle} which, in particular, rely on the root space decomposition from \cite[Sect. 2.5]{K}; see 
\cite[Lemmas 7.1, 8.3]{Gow} for more details.

\begin{lem}\label{lemma_last}
The center of ${\rm U}(\mathcal{L}(\A(n-1,n-1)))$ is trivial.
\end{lem}

Define the  {\em double Yangian $\dyp$ for $\A(n-1,n-1)$} as the quotient of $\dys$ over its $h$-adically completed ideal generated by the elements $b_{n|n}^{( \pm r)}$, where $r=1,2,\ldots.$ 
Due to \eqref{hto0maps},
the map in \eqref{isomorphism_sl} induces the isomorphism
$$
\dyp / h\dyp \cong  {\rm U}(\mathcal{L}(\A(n-1,n-1))). 
$$
Combining this isomorphism and Lemma \ref{lemma_last}, we obtain, analogously to Theorem \ref{centpro},

\begin{pro} 
The coefficients $b_{n|n}^{(\pm r)}$  with $r=1,2,\ldots$ are algebraically independent topological generators of the center of the double Yangian $\dys$.
Moreover, the center of  $\dyp$ is trivial.
\end{pro}

 \section*{Acknowledgment}
L.B. is member of Gruppo Nazionale per le Strutture Algebriche, Geometriche e le loro Applicazioni  (GNSAGA) of the Istituto Nazionale di Alta Matematica (INdAM).
This work has been supported  by Croatian Science Foundation under the project UIP-2019-04-8488. Furthermore, this work was supported by the project "Implementation of cutting-edge research and its application as part of the Scientific Center of Excellence for Quantum and Complex Systems, and Representations of Lie Algebras", PK.1.1.02, European Union, European Regional Development Fund. 

\linespread{1.0}


\begin{thebibliography}{9}
\bibitem{BaKo}
L. Bagnoli, S. Ko\v{z}i\'{c}, 
{\em Double Yangian and reflection algebras of the Lie superalgebra $\mathfrak{gl}_{m|n}$}, 
 Commun. Contemp. Math. (2024), https://doi.org/10.1142/S021919972450007X;
\href{https://arxiv.org/abs/2311.02410}{arXiv:2311.02410 [math.QA]}.


\bibitem{BK}
J. Brundan, A. Kleshchev, 
{\em Parabolic presentations of the Yangian $Y (\mathfrak{gl}_n )$}, 
Comm. Math. Phys.  \textbf{254} (2005), 191--220; 
\href{https://arxiv.org/abs/math/0407011}{arXiv:math/0407011 [math.QA]}.

\bibitem{CH}
H. Chang, H. Hu,
{\em A note on the center of the super Yangian $Y_{M|N}(\mathfrak{s})$},
J. Algebra \textbf{633} (2023), 648--665;
\href{https://arxiv.org/abs/2301.08100}{arXiv:2301.08100 [math.QA]}.

\bibitem{GGRW}
I. Gelfand, S. Gelfand, V. Retakh,  R. L. Wilson, 
{\em Quasideterminants},
 Adv. Math. \textbf{193} (2005), 56--141;
\href{https://arxiv.org/abs/math/0208146}{arXiv:math/0208146 [math.QA]}.

\bibitem{GR}
I. Gelfand, V. Retakh, 
{\em Quasideterminants, I}, 
 Selecta Math. (N.S.), \textbf{3} (1997), 517--546;
\href{https://arxiv.org/abs/q-alg/9705026}{arXiv:q-alg/9705026}.

\bibitem{Gow0}
 L. Gow, 
{\em On the Yangian $\text{Y}(\mathfrak{gl}_{m|n} )$ and its quantum Berezinian}
, Czech. J. Phys. \textbf{55} (2005), 1415--1420;
\href{https://arxiv.org/abs/math/0501041}{arXiv:math/0501041 [math.QA]}.

\bibitem{Gow}
L. Gow, 
{\em Gauss Decomposition of the Yangian $Y(\mathfrak{gl}_{m|n})$}, Comm. Math. Phys. \textbf{276} (2007), 799--825;
\href{https://arxiv.org/abs/math/0605219v3}{arXiv:math/0605219 [math.QA]}.

\bibitem{HM}
C. Huang, E. Mukhin,
{\em The duality of $\mathfrak{gl}_{m|n}$ and $\mathfrak{gl}_k$ Gaudin models},
J. Algebra \textbf{548} (2020), 1--24;
\href{https://arxiv.org/abs/1904.02753}{arXiv:1904.02753 [math.QA]}.

\bibitem{I}
K. Iohara,
{\em Bosonic representations of Yangian double $DY_{\hbar}(\mathfrak{g})$ with $\mathfrak{g}=\mathfrak{gl}_N,\mathfrak{sl}_N$},
J. Phys. A \textbf{29} (1996), 4593--4621;
\href{https://arxiv.org/abs/q-alg/9603033}{arXiv:q-alg/9603033}.

\bibitem{JY}
N. Jing, F. Yang,
{\em Center of the Yangian double in type A},
 Sci. China Math. (2024), https://doi.org/10.1007/s11425-022-2142-9;
\href{https://arxiv.org/abs/2207.01712v2}{arXiv:2207.01712 [math.QA]}.

\bibitem{K}
V. G. Kac, 
{\em Lie superalgebras}, 
Adv. Math. \textbf{26} (1977), 8--96.

\bibitem{Lu2209}
K. Lu,
{\em On Bethe eigenvectors and higher transfer matrices for supersymmetric spin chains},
J. High Energy Phys. \textbf{4} (2023),    120, 39 pp;
\href{https://arxiv.org/abs/2209.14416}{arXiv:2209.14416 [math-ph]}.

\bibitem{Lu2308}
K. Lu,
{\em Isomorphism between twisted $q$-Yangians and affine $\imath$quantum groups: type AI},
\href{https://arxiv.org/abs/2308.12484}{arXiv:2308.12484 [math.QA]}.

\bibitem{LM}
K. Lu, E. Mukhin, 
{\em Jacobi--Trudi identity and Drinfeld functor for super Yangian},
Int. Math. Res. Not. IMRN \textbf{2021} (2021), 16751–16810;
\href{https://arxiv.org/abs/2007.15573v2}{arXiv:2007.15573 [math.QA]}.

\bibitem{M_book}
A. Molev, 
{\em Yangians and classical Lie algebras}, 
Mathematical Surveys and Monographs, \textbf{143}. American Mathematical Society, Providence, RI, 2007.

\bibitem{MNO}
A. Molev, M. Nazarov, G. Olshanski, 
{\em Yangians and classical Lie algebras}, 
Russian Math. Surveys \textbf{51} (1996), 205--282;
\href{https://arxiv.org/abs/hep-th/9409025}{arXiv:hep-th/9409025}.

\bibitem{MR}
A. Molev, E. Ragoucy, 
{\em The MacMahon Master Theorem for right quantum superalgebras and higher Sugawara operators for $\widehat{gl}(m|n)$}, 
Mosc. Math. J. \textbf{14} (2014), 83--119;
\href{https://arxiv.org/abs/0911.3447}{arXiv:0911.3447 [math.RT]}.

\bibitem{Naz}
M. L. Nazarov, 
{\em Quantum Berezinian and the classical Capelli identity}, 
Lett. Math. Phys.
\textbf{21} (1991), 123--131.

\bibitem{NT}
M. Nazarov, V. Tarasov, 
{\em Representations of Yangians with Gelfand--Zetlin bases}, 
 J. Reine Angew. Math.  \textbf{496} (1998), 181--212;
\href{https://arxiv.org/abs/q-alg/9502008}{arXiv:q-alg/9502008}.

\bibitem{T}
 A. Tsymbaliuk, 
{\em Shuffle algebra realizations of type $A$ super Yangians and quantum affine
superalgebras for all Cartan data}, 
Lett. Math. Phys. \textbf{110} (2020), 2083--2111;
\href{https://arxiv.org/abs/1909.13732}{arXiv:1909.13732 [math.RT]}.

\bibitem{Z}
Y-Z. Zhang, 
{\em Super-Yangian double and its central extension},
Phys. Lett. A \textbf{234} (1997), 20--26;
\href{https://arxiv.org/abs/q-alg/9703027}{arXiv:q-alg/9703027}.
\end{thebibliography}
\end{document}